\input amstex
\documentstyle{amsppt}
\magnification=1200
\vcorrection{-0.3in}
\NoRunningHeads
\NoBlackBoxes
\topmatter
\title Positively curved manifolds with maximal discrete  
symmetry rank
\endtitle
\author Fuquan Fang \footnote{Supported by CNPq of Brazil, NSFC Grant
19741002, RFDP and Qiu-Shi Foundation of China
\hfill{$\,$}}\,\& \,Xiaochun Rong\footnote{Supported partially by
NSF Grant DMS 0203164 and a research found from Capital normal
university
\hfill{$\,$}}
\endauthor
\address Nankai Institute of Mathematics, Nankai University, Tianjing
300071, P.R.China.
\newline
Instituto de Matematica, Universidade Federal Fluminense, Niteroi, RJ 24005,
Brasil
\endaddress
\email ffang\@nankai.edu.cn  \hskip 4mm fuquan\@mat.uff.br
\endemail
\address Mathematics Department, Capital Normal University,
Beijing, P.R.C.
\newline .\hskip3mm Mathematics Department, Rutgers University, New
Brunswick,
NJ 08903, U.S.A
\endaddress
\email rong\@math.rutgers.edu
\endemail
\abstract Let $M$ be a closed simply connected $n$-manifold of positive 
sectional curvature. We determine its homeomorphism or homotopic type 
if $M$ also admits an isometric elementary $p$-group action 
of large rank. Our main results
are: There exists a constant $p(n)>0$ such that (1) If $M^{2n}$ admits an 
effective isometric $\Bbb Z_p^k$-action for a prime $p\ge p(n)$, then
$k\le n$ and ``$=$'' implies that $M^{2n}$ is homeomorphic to a 
sphere or a complex projective space. (2) If $M^{2n+1}$ admits an isometric 
$S^1\times \Bbb Z_p^k$-action for a prime $p\ge p(n)$, then 
$k\le n$ and ``$=$'' implies that $M$ is homeomorphic to a sphere.
(3) For $M$ in (1) or (2), if $n\ge 7$ and $k\ge 
[\frac{3n}4]+2$, then $M$ is homeomorphic to a sphere or homotopic to 
a complex projective space. 
\endabstract
\endtopmatter
\document

\vskip2mm

\head 0. Introduction
\endhead

\vskip4mm

A basic problem in Riemannian geometry is the classification of the positively curved manifolds whose isometry groups are large (cf. [Gro]). The largeness 
is measured by the dimension (resp. the rank of a maximal torus) of the
isometry 
group or the dimension of its orbit space. A typical example is the
classification 
of homogeneous manifolds of positive curvature ([AW], [Be], [Ber], [Wa]). 

Recently, there has been considerable progress on the classification of
positively 
curved manifolds with large symmetry rank (cf. [FMR], [FR1,2], [GS], [Hi],
[HK], 
[Ro2], [PS], [Wi], [Ya], etc). The symmetry rank of a Riemannian manifold is
the 
rank of a maximal torus of the isometry group. 

The maximal  rank theorem of Grove-Searle asserts that if a closed  
positively curved $n$-manifold $M$ admits an isometric torus $T^k$-action,
then 
$k\le [\frac{n+1}2]$ (the integer part) and ``$=$'' if and only if $M$ is 
diffeomorphic to a sphere, a lens space or a complex projective space ([GS]). Note that any $G$-action considered in this paper is assumed to be {\it
effective}.

A natural next step is to (homeomorphically) classify closed simply connected positively curved manifolds of the {\it almost maximal} 
symmetry rank $[\frac{n+1}2]-1=[\frac{n-1}2]$. This was previously done for $n=4$ by Hsiang-Kleiner ([HK]), for $n=5$ by Rong ([Ro2]). Recently, 
the homeomorphism classification for all $n\ge 8$ was obtained by 
Fang-Rong ([FR1], cf. [FMR], [Wi]). And these manifolds are 
{\it homeomorphic} to a sphere or a complex projective space 
or a quaternionic projective plane (for $n=8$). 
The situation in the remaining dimensions $n=6, 7$ are quite 
subtle and complicated; there are many examples which are not 
homeomorphic to any rank one symmetric space ([AW], [Es]). The homeomorphism classification in dimensions $6$ and $7$ has been under investigation ([FR1,2]).

Recently, Burkhard Wilking obtained the following almost $1/2$-maximal 
rank theorem: For $n\ge 10$, if  a closed simply connected positively curved $n$-manifold $M$ with symmetry rank at least 
$(n/4+1)$, then $M$ is homeomorphic to a sphere or a quaternionic projective space or {\it homotopy} equivalence to a complex projective space ([Wi]).

A new tool used in [Wi] (also in [FR]) is a connectedness theorem 
of Wilking: If $N$ is a closed totally geodesic $r$-submanifold of 
an $n$-manifold $M$ of positive sectional curvature, then the 
homotopy group $\pi_i(M,N)=0$ for $i\le 2r-n+1$ (see Theorem 1.11). Indeed, this theorem, together with other connectedness theorems 
(see Theorems 1.9-1.12) will play a crucial role in the present 
paper.

Hopefully, the above classification results together would cover 
almost all known examples of positively curved manifolds with 
large symmetry rank, except in dimension $13$ ([Ba1]). This should 
be helpful in looking for new examples of positively curved manifolds
(cf. [Gro]).

The purpose of this paper is to study the following more 
general problem (compare to Problem 24 in [Pe]): 
$$\text{\it Classify positively curved manifolds which admit 
isometric $\Bbb Z_p^k$-actions ($k$ large).}$$
One hopes that this will lead to new insights on positively curved 
manifolds with abelian symmetry (see Problems 0.1 and 0.2).

A guideline is to extend theorems on positively curved 
manifolds with (large) symmetry rank to positively curved manifolds 
with (large) {\it discrete} symmetry rank. The discrete $p$-symmetry
rank ($p$ is a prime) is the largest number $r$ such that the isometry group contains an elementary $p$-group of rank $r$. This is partly justified 
by a fascinating fact in the theory of 
compact transformation groups: a $T^k$-action shares many properties 
with an elementary $p$-group $\Bbb Z_p^k$-action 
(e.g. the Smith theorem and the Borel theorem) while these 
properties are generally false for any other compact Lie group (cf. [Hs]).

We begin to state the main results in this paper. The first two 
results (Theorems A and B) may be considered as an analog of the maximal  
rank theorem of Grove-Searle for elementary $p$-groups.

\vskip2mm

\proclaim{Theorem A (Maximal Discrete Rank)} 

There exists a constant $p(n)>0$ such that if a 
closed simply connected $2n$-manifold of positive sectional curvature
admits an isometric $\Bbb Z_p^k$-action with prime $p\ge p(n)$, then 
$k\le n$ and ``$=$'' implies that $M$ is homeomorphic 
to a sphere or a complex projective space.
\endproclaim

\vskip2mm

The number $p(n)$ can be chosen to be the Gromov's Betti number 
bound for closed  $n$-manifolds of non-negative sectional curvature 
(see Theorem 1.7).

A basic ingredient in the proof is that the $\Bbb Z^k_p$-fixed point set 
is not empty; see Lemma 2.1 (note that this may not hold when 
$p$ is small; the unit sphere $S^ 2$ admits isometric $\Bbb Z_2^3$-action by
reflections without fixed point). By comparing
this to the Isotropy rank theorem (cf. [Gro]), a natural problem arises:
 
\vskip2mm

\example{Problem 0.1} Let $M$ be a closed $(2n+1)$-manifold of 
positive sectional curvature. Show that for $p\ge p(n)$, any 
isometric $\Bbb Z_p^k$-action ($k\ge 2$) has a subgroup isomorphic to 
$\Bbb Z_p^{k-1}$ whose fixed point set is not empty.  
\endexample

\vskip2mm

Note that Problem 0.1 implies that the fundamental group of any positively 
curved $n$-manifold cannot be isomorphic to $\Bbb Z_p\oplus \Bbb Z_p$ for 
$p\ge p(n)$ (compare to [Ba2], [GSh], [Sh]). Indeed, this is 
also a consequence of the Almost cyclicity conjecture (cf. [Ro1]).

If one solves Problem 0.1, then following the argument in this 
paper one can prove a Maximal Discrete Rank theorem in odd-dimensions. 
We will prove the case of Problem 0.1 when $\Bbb Z_p^{k+1}$ is a subgroup 
of $S^1\times \Bbb Z_p^k$ (see Lemma 3.1). Hence, we are able to prove 

\vskip2mm 

\proclaim{Theorem B} 

Let $M$ be a closed simply connected $(2n+1)$-manifold of positive sectional curvature. If $M$ admits an isometric $S^1\times \Bbb Z_p^k$-action 
($p\ge p(n)$), then $k\le n$ and ``$=$'' implies that $M$ is 
homeomorphic to a sphere.
\endproclaim

\vskip2mm

Recall that the maximal  rank theorem of Grove-Searle follows easily 
from the theorem in [GS]: If a closed positively curved manifold $M$ admits an isometric circle action with fixed point set of codimension $2$, then $M$ 
is diffeomorphic to a sphere, a lens space, or a complex projective space. 
It is easy to see, using [Sm] and the connectedness theorem of Wilking 
mentioned earlier, that in the above, if one replaces the isometric 
$S^1$-action by any isometric $\Bbb Z_p$-action, then the universal 
covering of $M$ is homeomorphic to a sphere, provided $\dim(M)$ is odd 
(see Lemma 4.2).

\vskip2mm

\example{Problem 0.2} Let $M$ be closed simply connected $2n$-manifold of 
positive sectional curvature. Assume $M$ admits an isometric 
$\Bbb Z_p$-action ($p\ge p(n)$) with a fixed point set of codimension $2$. 
Show that $M$ is homotopic to a sphere or a complex projective space.
\endexample

\vskip2mm

Next, we consider the problem of a possible analog of Wilking's 
almost $1/2$-maximal rank theorem for the elementary $p$-groups. This 
is quite subtle due to that the almost $1/2$-maximal discrete symmetry 
rank seems to be an inadequate condition for the existence of 
a totally geodesic submanifold of codimension at most $\dim(M)/4$ 
(cf. [Wi]). Nevertheless, a discrete analog of Wilking's almost $1/2$-maximal 
rank theorem would be a $c$-maximal discrete rank theorem for some 
constant $\frac 12\le c< 1$. Our last two theorems imply 
that $c$ is not greater than $3/4$. 

\vskip2mm

\proclaim{Theorem C (Almost $3/4$-maximal Discrete Rank)}

For $n\ge 7$, if a closed simply connected $2n$-manifold of positive 
sectional curvature satisfies 
$$\text{Symrank}_p(M)\ge \cases [\frac{3n}4]+2 & n=1, 2 \text{ mod } 4\\
[\frac{3n}4]+1& \text{otherwise}\endcases\qquad \qquad p\ge p(n),$$
then $M$ is homeomorphic to a sphere or homotopic equivalent to 
a complex projective space.
\endproclaim

\vskip2mm

The following is a weak version of Theorem C in odd-dimensions; which 
also partially strengthens Theorem B.

\vskip2mm

\proclaim{Theorem D}

For $n\ge 7$, a closed simply connected $(2n+1)$-manifold of positive 
sectional curvature is homeomorphic to a sphere, if $M$ admits an isometric 
$S^1\times \Bbb Z_p^k$-action ($p\ge p(n)$) such that $k
\ge \cases [\frac{3n}4]+2& n=1\, \text{ mod } 4\\
[\frac{3n}4]+1 & \text{otherwise}\endcases$.
\endproclaim

\vskip2mm

\remark{Remark \rm 0.3} Note that in cases of Theorems C and D 
where there is an isometric $T^k$-action (resp. $T^{k+1}$-action),
the lower bound for $k$ can be further lowered by one, using the 
equivariant version of Theorem 1.11 in [Wi].
\endremark

\vskip2mm

We now give an indication about the proof of Theorem A (the proofs 
of Theorems C and D are in a similar spirit but are more involved). 

The proof is based on the following deep topological results: 
For $n\ge 4$ the generalized Poincare conjecture holds true
($n=4$ by Freedman and $n\ge 5$ by Smale, [Fr], [Sm]) and the
theory of Sullivan's characteristic varieties [Su] (cf. Theorem 1.3),
which implies that a homotopy complex projective space $M$
is homeomorphic to $\Bbb CP^n$ if $M$ contains a closed submanifold 
of codimension $2$ dual to a generator of $H^2 (M;\Bbb Z)\cong \Bbb Z$
which is homeomorphic to $\Bbb CP^{n-1}$ (see Lemma 2.4). 

Consider $M$ as in Theorem A. By employing the Borel theorem (cf. [Hs])
and the Gromov's Betti number estimate ([Gr]), we show that for $p\ge p(n)$, the $\Bbb Z^k_p$-fixed point set is not empty (Lemma 2.1). From the isotropy representation at a fixed point, it follows that $k\le n$,
and ``$=$'' implies that there is a sequence of closed totally geodesic
submanifolds, 
$$M^4\subset M^6\subset \cdots \subset M^{2n}=M,$$
such that $M^{2i}$ admits an isometric $\Bbb Z_p^i$-action. By the theorem of Wilking mentioned earlier (Theorem 1.11), it is easy to see that if 
$M^4$ is homotopic to a sphere or a complex projective space, 
then $M$ is homotopic to a sphere or a complex projective space. 
Indeed, if $M^4$ is homeomorphic to a sphere or complex projective 
space, then by the above topological results one upgrades the homotopy 
equivalence to that $M$ is homeomorphic to a sphere or a complex projective 
space. We then complete the proof by showing that $M^4$ is homeomorphic to a sphere or a complex projective space (Proposition 2.3).

The rest of the paper is organized as follows:

\noindent In Section 1, we collect results that are used in the proofs of 
Theorems A-D.

\noindent In Section 2, we prove Theorem A.

\noindent In Section 3, we prove Theorem B.

\noindent In Section 4, we prove Theorems C and D.

\vskip4mm

{\bf Acknowledgment}: The authors would like to thank Burkhard Wilking 
for some useful comments on the early version 
and thank Wolfgang Ziller for some references. The second author 
would like to thank Burkhard Wilking for informing him of 
the main results in [Wi].

\vskip10mm

\head 1. Preliminaries
\endhead

\vskip4mm

In this section, we will collect theorems that will be used in the 
rest of the paper.

\vskip4mm

\subhead a. Some deep topological theorems
\endsubhead

\vskip4mm

The generalized Poincar\'e conjecture says that any homotopy $n$-sphere is
a homeomorphic sphere. By [Fr] and [Sm], this conjecture remains open 
only for $n=3$.

\vskip2mm

\proclaim{Theorem 1.1 (Freedman)} Every homotopy $4$-sphere is
homeomorphic to a sphere.
\endproclaim

\vskip2mm

\proclaim{Theorem 1.2 (Smale)} For $n\ge 5$, every homotopy $n$-sphere is
homeomorphic to a sphere.
\endproclaim

\vskip2mm

In general, the homotopy class of a closed manifold may contain
many homeomorphism types. In [Su], Sullivan classified the homeomorphic
types of a homotopy complex projective space.

Let $M$ be a manifold with the homotopy type of $\Bbb CP^n$. Let 
$h: M\to \Bbb CP^n$ be a homotopy equivalence. For any $\Bbb CP^i
\subset \Bbb CP^n$, let $h^{-1}(\Bbb CP^i)\subset M$ denote the 
transverse submanifold (the preimage of a map homotopic to $h$ 
and transverse to $\Bbb CP^i$.) Let $\sigma _h(i)$ be the signature 
(an integer) of $h^{-1}(\Bbb CP^i)$ if $i$ is even (an integer), 
and the Kervaire invariant (a mod $2$ integer) if $i$ is odd 
(cf. [Su]). The invariant $\sigma _h(i)$ depends only on the 
homotopy class $h$. Indeed, if $N^{2i}\subset M$ is another 
manifold homologous to $h^{-1}(\Bbb CP^i)$, then $\sigma _h(i)$ 
is equal to the signature (resp. Kervaire invariant) of $N^{2i}$. 

\vskip2mm 

\proclaim{Theorem 1.3 (Sullivan)} Let $hS(\Bbb CP^n)$ denote the set of 
closed manifolds homotopy equivalent to $\Bbb CP^n$. For $i=1, \cdots,
[n/2]$,
$\sigma _h(i)$ defines an one-to-one correspondence between $hS(\Bbb CP^n)$ 
and the set $\prod_{i=1}^{[n/4]} (\Bbb Z\times\Bbb Z_2)$.
\endproclaim

\vskip2mm

Recall that the Kervaire invariant is zero if a manifold of 
dimension $(4i+2)$ has no middle dimensional homology, e.g. 
manifolds homotopy equivalent to $\Bbb C P^{2i+1}$. On the 
other hand, the signature is a homotopy invariant.  As an 
immediate corollary of Theorem 1.3, one gets 

\vskip2mm

\proclaim{Corollary 1.4} If $M\in hS(\Bbb CP^n)$ has a codimension 
$2$ submanifold homeomorphic to $\Bbb CP^{n-1}$ which represents a 
generator of $H_{2n-2}(M)$, then $M$ is homeomorphic to $\Bbb CP^n$.
\endproclaim

\vskip4mm

\subhead b. Fixed point sets
\endsubhead

\vskip4mm

Let $G$ denote a compact Lie group. For a $G$-space, let $F(G,M)$ denote 
the fixed point set. Given any $G$-invariant metric, each component 
of $F(G,M)$ is a totally geodesic submanifold.  

It is well known that for $G=T^k$ or $\Bbb Z_p^k$ ($p$ is a prime), the 
topology of $F(G,M)$ is closely related the topology of $M$. For 
instance, the Euler characteristic $\chi(M)=\chi(F(T^k,M))$. In the proof 
of Theorem A, the following results will be used (cf. [Hs]).

\vskip2mm

\proclaim{Lemma 1.5} Let $M$ be a closed  $\Bbb Z_p$-space ($p$ is a prime).

Then 
$$\chi(M)=\chi(F(\Bbb Z_p,M))\,(\text{mod}\,\, p).$$
\endproclaim

\vskip2mm

\proclaim{Theorem 1.6 (Borel)} Let $G=T^k$ or $\Bbb Z_p^k$ ($p$ is a prime), and let $M$ be a $G$-space. Then 
$$\text{rank}(H_*(F(G,M),\ell))\le \text{rank}(H_*(M,\ell)),$$
where $\ell=\Bbb R$ if $G=T^k$ and $\Bbb Z_p$ otherwise. 
\endproclaim

\vskip2mm

By Theorem 1.6, the number of components of $F(G,M)$ is bounded above by 
the total Betti number of $M$. If $M$ also has non-negative sectional 
curvature, then Gromov Betti number estimate ([Gr]) implies that 
the number of components of $F(G,M)$ is bounded above a constant 
depending only on $\dim(M)$. 

\vskip2mm

\proclaim{Theorem 1.7 (Gromov)} Let $M$ be a closed  $n$-manifold of 
non-negative sectional curvature. Then for any coefficient field $\ell$, 
the total Betti number, $\text{rank}(H_*(M,\ell))\le b(n)$.
\endproclaim

\vskip4mm

\subhead c. The connectedness principle of positive curvature
\endsubhead

\vskip4mm

As seen in the introduction, in the study of positively curved 
manifolds with large isometry groups, the recent break-through 
is the Wilking's connectedness theorem (see Theorem 1.11). 
His proof is a combination of the Morse theory and 
the standard Synge type argument. Soon after Wilking announced 
the results of [Wi], a connectedness principle of positive 
curvature was obtained in [FMR] (see Theorem 1.8) 
which not only gives a uniform formulation for the Synge 
theorem (Theorem 1.12), the Frankel theorem (Theorem 1.9) 
and the Wilking theorem (Theorem 1.11), but also includes 
a new regularity theorem (which is not 
required in this paper). The connectedness principle of 
positive curvature can be viewed as an analog of the 
connectedness principle in the algebraic geometry, 
for details, see [FMR] and references within.  

An immersion is called {\it totally geodesic} if the second 
fundamental form vanishes.

\vskip2mm 

\proclaim{Theorem 1.8 (Connectedness Principle of positive 
curvature)} Let $M$ be a closed $m$-manifold of positive sectional curvature,
and
let $\Delta$ denote the diagonal of $M\times M$.
Assume $N=N_1\times N_2$ and $f=(f_1,f_2): N\to M\times M$, 
where $f_i: N_i\to M$ is totally geodesic $n_i$-dimensional 
immersion. Then $(n=n_1+n_2$)
  
\noindent (1.8.1) If $n\ge m$, then $f^{-1}(\Delta)$ is non-empty.  
  
\noindent (1.8.2) If $n\ge m+1$ and $M$ is simply connected, then  
$f^{-1}(\Delta )$ is connected.  
  
\noindent (1.8.3) For $n\ge m+i$ there is an exact sequence  
$$\pi _i(f^{-1}(\Delta ))\to \pi _i(N) @>(p_1f)_* - (p_2f)_*>>  
\pi _i(M)\to \pi _{i-1}(f^{-1}(\Delta ))\to \cdots .$$  
  
\noindent (1.8.4) There are natural isomorphisms,  
$\pi_i(N_1,f^{-1}(\Delta ))\to\pi_i(M,N_2)$ and 
$\pi_i(N_2,f^{-1}(\Delta))\to \pi_i(M,N_1)$ 
for $i\le n-m$ and surjections for $i=n-m+1$,   
where $\pi_i(N_j,f^{-1}(\Delta))$ is understood as the $i$-th  
homotopy group of the composition map  
$f^{-1}(\Delta)\subset  N@>p_j>> N_j$.  
\endproclaim  

\vskip2mm

Note that (1.8.1) is a strengthened version of the Frankel theorem
(stated below) where embedded submanifolds are assumed. 

\vskip2mm

\proclaim{Theorem 1.9 (Frankel)} Let $M$ be a closed  manifold of positive 
sectional curvature, and let $f_i: N_i\to M$ be two closed totally geodesic 
immersions. If $\dim(N_1)+\dim(N_2)\ge \dim(M)$, then $f_1(N_1)\cap f_2(N_2)
\ne \emptyset$.
\endproclaim

\vskip2mm

Note that (1.8.4) immediately implies the following:

\vskip 2mm

\proclaim{Theorem 1.10} Let $M$ be a closed $m$-manifold of 
positive curvature. Let $N_1, N_2$ be embedded totally geodesic 
submanifolds in $M$ of dimensions $n_1, n_2$ respectively. 
Set $n=n_1+n_2$. Then there are natural isomorphisms
$$\pi_i(N_1, N_1\cap N_2)\to \pi_i(M,N_2),\qquad \qquad
\pi_i(N_2, N_1\cap N_2)\to \pi_i(M,N_1)$$
for $i\le n-m$, and surjections for $i=n-m+1$.
\endproclaim

\vskip2mm

A map from $N$ to $M$ is called {\it $(i+1)$-connected}, if it induces
an isomorphism up to the $i$-th homotopy group and a surjective
homomorphism on the $(i+1)$-th homotopy group.

Theorem 1.10 implies the following result of Wilking ([FMR]). 

\vskip2mm

\proclaim{Theorem 1.11 (Wilking)} Let $M$ be a closed $m$-manifold 
of positive sectional curvature. Let $N_1$ and $N_2$ be two totally 
geodesic submanifolds of dimensions $n_1$ and $n_2$. If $n_2\ge n_1$, 
then the inclusions,

\noindent (1.11.1) $i_2: N_2\hookrightarrow M$ is $(2n_2-m+1)$-connected.

\noindent (1.11.2) $N_1\cap N_2\hookrightarrow N_1$ is 
$(n_1+n_2-m)$-connected.
\endproclaim

\vskip2mm

Indeed, the classical Synge theorem can also be formulated as a
connectedness
theorem (compare (1.8.1)). This was pointed to us by Karsten Grove.

\vskip2mm

\proclaim{Theorem 1.12 (Synge)} Let $M$ be a closed orientable $n$-manifold 
of positive sectional curvature, and let $\phi$ be an isometry. Let $f: 
M\to M\times M$ be a totally geodesic embedding by $f(x)=(x,\phi(x))$. 
Then $f^{-1}(\Delta)$ is not empty under the following situations:

\noindent (1.12.1) $n$ is even and $\phi$ is orientation preserving.

\noindent (1.12.2) $n$ is odd and $\phi$ is orientation reversing.
\endproclaim

\vskip2mm

There are other connectedness theorems that follow from Theorem 1.8 
(see [FMR]); we do not state them here since they are not required 
in this paper.

\vskip4mm

\subhead d. Alexandrov spaces with positive curvature
\endsubhead

\vskip4mm

Recall that an Alexandrov space, $X$, is a finite Hausdorff 
dimensional complete metric space with a lower curvature 
bound in distance comparison sense, cf. [BGP]. In particular, a 
Riemannian manifold of sectional curvature bounded from below is 
an Alexandrov space. For an Alexandrov space, one can define 
orientability in a standard way using atlas of distance maps; 
see [Pet].

The following version of Theorem 1.12 for Alexandrov spaces of 
positive curvature will be used in the proof of Theorem B. 

\vskip2mm 

\proclaim{Theorem 1.13 (Petrunin)} Let $X$ be an orientable Alexandrov 
space of cur$(X)\ge 1$. If dim$(X)$ is even, then any 
orientation-preserving isometry has a fixed point. In particular, 
$X$ is simply connected.
\endproclaim

\vskip2mm

Recall that the quotient space of a Riemannian manifold by an 
isometric group action is, with the quotient metric, not 
necessarily a Riemannian manifold (it may not even be a manifold). 
In the comparison, an advantage with an Alexandrov space is 
(cf. [BGP]):   

\vskip2mm

\proclaim {Lemma 1.14} Let $X$ be an Alexandrov space with curvature 
$\ge -\Lambda$. Let $G$ be a compact group of isometries. Then, the 
quotient space, $X/G$, is also an Alexandrov space with curvature 
$\ge-\Lambda$.
\endproclaim

\vskip4mm

\head 2. Proof of Theorem A
\endhead

\vskip4mm

For an isometric $T^k$-action on a closed  positively 
curved $n$-manifold, a basic property is the Isotropy 
rank theorem: 
If $n$ is even, then the fixed point set is not empty, 
and if $n$ is odd, then there is a circle orbit (cf. [Ko], [Ro2]). 
This also plays a role in the proof of the maximal rank theorem 
of Grove-Searle ([GS]). 

The following result may be considered as a discrete version of
the Isotropy rank theorem in even-dimensions (compare to Lemma 3.1). 

\vskip2mm  

\proclaim{Lemma 2.1 (Discrete Isotropy Rank)} Assume a closed  
positively curved $2n$-manifold $M$ admits an isometric $\Bbb 
Z_p^k$-action with $p\ge p(n)$. Then the fixed point set 
is not empty and $k\le n$ and `$=$' implies that the 
fixed point set is finite.
\endproclaim

\vskip2mm

Note that the requirement for a lower bound on $p$ is of necessary; for 
instance the unit sphere $S^{2n}$ admits the obvious isometric 
$\Bbb Z_2^{2n+1}$-action without fixed point.

\vskip2mm

\demo{Proof} Without loss of the generality, we may assume 
that $p(n)>2$. We shall determine the value of $p(n)$ at the end of 
the proof.

We first assume that $M$ is orientable. Let $H_p$ denote any 
$\Bbb Z_p$ subgroup of $\Bbb Z_p^k$ generated by $\alpha$. Note that
$\alpha^2$ is 
always orientation-preserving and $F(H_p,M)=F(\alpha^2,M)$ since
$H_p$ is generated by $\alpha^2$.
By Theorem 1.12, we can assume an $H_p$-fixed point component $F_0\ne 
\emptyset$. Note that $F_0$ is a closed totally geodesic submanifold
of even-codimension. We claim that $F_0$ is preserved by $\Bbb Z_p^k$.
Then, by induction on $n$ one can conclude that $\Bbb Z_p^k$ has a fixed point in $F_0$. If not, $F(H_p,M)$ must have at least $p$ components. 
By Theorem 1.6 and Theorem 1.7, we then obtain  
$$\split p&\le \#\{\text{components of $F(\Bbb Z_p,M)$}\}\\&\le \text{rank}
(H_*(F(\Bbb Z_p,M),\Bbb Z_p))\\& \le \text{rank}(H_*(M,\Bbb Z_p))\le
b(n).\endsplit$$
By now it is clear that if one chooses $p(n)=b(n)+1$, we then see 
a contradiction.

Let $F$ be a $\Bbb Z^k_p$-fixed point component. Then $\dim (F)=2(n-r)$
where $r$ is a positive integer. Then, $\Bbb Z_p^k$ acts on the normal space
of $F$ at $x$ by isometries.
We can consider $\Bbb Z_p^k$ as a subgroup of the maximal torus of
$O(2r)$ and therefore $k\le r\le n$ and $k=n$ implies that $r=n$ and 
thus $\dim (F)=0$.

If $M$ is not orientable, then the double covering $\tilde M$ of $M$ is
simply
connected (Theorem 1.12). Hence, we can apply the above to the lifting
$\Bbb Z_p^k$ action on $\tilde M$ and conclude the desired result.
\qed\enddemo

\vskip2mm

The maximality has the following property.

\vskip2mm

\proclaim{Lemma 2.2} Let the assumptions be as in Theorem A with $k=n$. 
Then there are closed totally geodesic submanifolds, $M^4\subset M^6
\subset\cdots \subset M^{2n}=M$ such that $M^{2i}$ admits an isometric 
$\Bbb Z_p^i$-action. 
\endproclaim

\vskip2mm

\demo{Proof} 
Consider the linear isotropy representation of $\Bbb Z_p^n$ at a 
fixed point $x\in M$. Let $\Bbb Z_p^{n-i} \subset \Bbb Z_p^n$ be the product
subgroup of the first $(n-i)$-factors. Let $M^{m_i}$ be the connected 
fixed point submanifold of $\Bbb Z_p^{n-i}$ containing $x$. It suffices to
verify the dimension(s) of the fixed point manifolds are exactly $2i$.

Note that the tangent space of $M^{m_i}$ at $x$ is  the invariant linear
subspace 
of the $\Bbb Z_p^{n-i}$-action on $T_xM$.  
By the proof of Lemma 2.1 we know that $m_i\le 2i$.  Consider the action of
$\Bbb Z_p^n$ on 
$F_i$. This action has the principal isotropy group $\Bbb Z_p^{n-i}$.
Therefore
the induced action by the quotient group $\Bbb Z_p^i=\Bbb Z_p^n/ \Bbb
Z_p^{n-i}$ is 
effective. This implies that $\Bbb Z_p^i$ embeds into the orthogonal group
$O(m_i)$ as
a subgroup. Therefore $m_i\ge 2i$. The desired result follows.
\qed\enddemo

\vskip2mm

We first prove Theorem A for $n=4$. For a closed  manifold $M$, 
let $\chi(M)$ denote the Euler characteristic of $M$.

\vskip2mm

\proclaim{Proposition 2.3} Let $M^4$ be a closed simply connected 
manifold of positive sectional curvature. If $M^4$ admits an isometric $\Bbb Z_p^2$-action with $p>\chi(M^4)$, then $M^4$ is homeomorphic to $S^4$ or $\Bbb CP^2$. 
\endproclaim

\vskip2mm

Note that one can alway apply Proposition 2.3 if $p>b(4)$, the 
constant from Theorem 1.7 for $n=4$. 

\vskip2mm

\demo{Proof of Proposition 2.3} We first show that the fixed point set $F(\Bbb
Z_p^2,M^4)
\ne \emptyset$. Since $F(\Bbb Z_p^2,M^4)$ has even codimension, from 
the isotropy representation $F(\Bbb Z_p^2,M^4)$ must be a finite set.

We will argue by contradiction. By Theorem 1.12, we can assume 
a $\Bbb Z_p$-subgroup, $H_p$, of $\Bbb Z_p^2$ such that $F(H_p,M^4)
\ne \emptyset$. If $\dim(F(H_p,M^4))=2$, then $F(H_p,M^4)$ has only 
one component of dimension $2$ and therefore $\Bbb Z_p^2$ must 
preserve this component (see Theorem 1.9). Since the component
is a closed totally geodesic submanifold (which is homeomorphic to 
$S^2$), it contains $\Bbb Z^2_p$-fixed points; a contradiction. 
If $F(H_p,M^4)$ is a finite set, since $F(\Bbb Z_p^2,M^4)
=\emptyset$, $|F(H_p,M^4)|\ge p$. By Theorem 1.6, we obtain
$$p\le \text{rank}(H_*(F(H_p,M^4),\Bbb Z_p)\le \text{rank}
(H_*(M^4,\Bbb Z_p))=\chi(M^4)<p,$$
a contradiction.

We then show that $\chi(M^4)\le 3$. By [Fr], it follows that $M^4$ is 
homeomorphic to $S^4$ or $\Bbb CP^2$. 

Observe that at each $\Bbb Z_p^2$-fixed point, from the isotropy 
representation we see that there are two $\Bbb Z_p$-isotropy groups 
whose fixed point sets are two spheres and their intersection is 
this $\Bbb Z_p^2$-fixed point. Hence, we can use a one-dimensional
complex to represent the singular set in the orbit space $M^4/\Bbb Z_p^2$,
where every edge indicates a two sphere fixed point set of some
isotropy group. The $\Bbb Z_p^2$-fixed point are the vertices and at each
vertex there are two edges. If there are more than three vertices, 
then there must exist two edges which do not share any vertex. This implies
that there are two totally geodesic $2$-spheres in $M^4$ which 
do not intersect, a contradiction (see Theorem 1.9). On the 
other hand, the number of vertices cannot be one. 

By Lemma 1.5, the above implies that $\chi(M^4)=2+kp<p$ or 
$3+kp<p$ with $k\ge 0$ since 
$\chi(M^4)=2+b_2\ge 2$. This implies that 
$\chi(M^4)=2$ or $3$.
\qed\enddemo

\vskip2mm 

\proclaim {Lemma 2.4} Let $M$ be a homotopy complex projective
space. If $M$ has a codimension $2$ submanifold $N$ homeomorphic 
to $\Bbb CP^{n-1}$ such that the inclusion map is a $3$-connected, 
then $N$ represents a generator of $H_{2n-2}(M)$ and therefore 
$M$ is homeomorphic to $\Bbb CP^n$.
\endproclaim

\vskip2mm

\demo{Proof} By Corollary 1.4, it suffices to show that $N$ represents
a generator of $H_{2n-2}(M)$. Let $x\in H^2(M)$ denote a generator. Recall 
that the 
cohomology ring $H^*(M)\cong H^*(\Bbb CP^n)\cong \Bbb Z[x]/(x^{n+1}=0)$. 
Since $i$ is a $3$-equivalence, $i: N\to M$ induces an isomorphism 
$i^*: H^2(M)\to H^2(N)$. If $i_*([N]) \in H_{2n-2}(M)$
represents $d$ times the generator, then $i^* (x^{n-1})[N]=d$. On the
other  hand, since the cohomology ring  $H^*(N)=\Bbb Z[i^*(x)]/
(i^* (x)^{n}=0)$. Thus $i^* (x^{n-1})[N]=\pm 1$. This proves that $d=\pm
1$, i.e. $i _*([N])$ is a generator of $H_{2(n-1)}(M)$.
\qed\enddemo

\vskip2mm

With the above preparation, we are ready to prove Theorem A.

\vskip2mm

\demo{Proof of Theorem A} 

We first show that $M$ is homotopy equivalent to a sphere or a complex 
projective space.

By Lemma 2.1, we obtain $k\le n$. By Lemma 2.2, we obtain a sequence 
of closed orientable totally geodesic submanifolds,
$$M^4\subset M^ 6\subset \cdots \subset M^{2n}=M.$$
By Theorem 1.11, we orderly conclude that $M^{2n-2},...,M^4$ are simply
connected.
By Proposition 2.3, we conclude that $M^4$ is homeomorphic to 
$S^4$ or $\Bbb CP^2$. By Theorem 1.11 again, we conclude 
that $M^6$ satisfies that $\pi _2(M^6)=\Bbb Z$ (resp. $0$) if $M^4=\Bbb CP^2$(resp. $S^4$) and $\pi_3(M)=0$. By the Poincar\'e duality and the Hurewicz 
theorem, we conclude that $M^6$ is homotopy equivalent 
to $\Bbb CP^3$ (resp. $S^6$) if $M^4\overset{\text{homeo}}\to \simeq 
\Bbb CP^2$ (resp. $S^4$). 
By iterating the same argument to $M^8,..., M^{2n}$, 
we conclude that $M^{2n}$ is homotopic to a sphere or a 
complex projective space.

Finally, if $M$ is a homotopy sphere, then by Theorems 1.1 and 1.2 $M$ 
is a homeomorphic sphere for $n\ge 4$. 

If $M$ is a homotopy complex projective space, by Proposition 2.3, Lemma 2.4 and Theorem 1.11 we then conclude that $M^6$ is homeomorphic to a complex 
projective space. Iterating this $(n-1)$-times, we eventually conclude that
$M^{2n}=M$ 
is homeomorphic to a complex projective space. 
\qed\enddemo

\vskip4mm

\head 3. Proof of Theorem B
\endhead

\vskip4mm

The proof of Theorem B follows the same strategy as in the proof of
Theorem A.

The first lemma may be considered as a weak version of Lemma 2.1 in
odd-dimensions. 

\vskip2mm

\proclaim{Lemma 3.1} Let $M$ be a closed simply connected 
$(2n+1)$-manifold of positive sectional curvature. If $M$ 
admits an isometric $S^1\times \Bbb Z_p^k$-action ($p\ge p(n)$)
then $k\le n$ and $\Bbb Z_p^{k+1}\subset S^1\times \Bbb Z_p^k$ 
has a $\Bbb Z^k_p$-subgroup whose fixed point set is not empty
and contains a circle orbit, and `$=$' implies that the 
fixed point set has dimension one.
\endproclaim

\vskip2mm

\proclaim{Sublemma 3.2} Let the assumptions be as in Lemma 3.1.
Then the $S^1\times \Bbb Z_p^{k}$-action has a finite isotropy 
group.
\endproclaim

\vskip2mm

\demo{Proof} Let $p: M\to M^*:=M/S^1$ denote the orbit
projection, where $S^1=S^1\times \{\text{id}\}\subset 
S^1\times \Bbb Z_p^k$.
Since $M$ is orientable and since the $S^1$-action preserves 
the orientation, $M^*$ is an orientable
Alexandrov space of positive curvature (Lemma 1.14). Note that
the $\Bbb Z_p^k$-action on $M$ descends to an effective
$\Bbb Z_p^k$-action on $M^*$.

Let $H\cong \Bbb Z_p$ be any subgroup of $\Bbb Z_p^k$.
By Theorem 1.13 we can assume that $H$ has a fixed point
$x^*\in M^*$. Let $S^1(x)$ denote the $S^1$-orbit at $x$
such that $p(x)=x^*$. Then $S^1\times H$ preserves $S^1(x)$
(which is either a circle or a point). This implies that the
isotropy group $H_x$ (of the $S^1\times \Bbb Z_p^k$-action) at $x$
contains a subgroup isomorphic to $H$.

If $S^1(x)$ is a circle, then the isotropy group of $S^1\times
\Bbb Z_p^k$ at $x$ is finite. Otherwise, $H_x$ contains
$S^1\times H$. Consider the isotropy representation of $H_x$
at the tangent space at $x$. From standard linear algebra,
one concludes that $H_x$ has a finite isotropy group on the
tangent space.
\qed\enddemo

\vskip2mm

\demo{Proof of Lemma 3.1}

We proceed by induction on $n$, $\dim(M)=2n+1$. Note that Lemma 3.1
holds
trivially for $n=1$. Assume Lemma 3.1 holds for $\dim(M)<2n+1$.

By Sublemma 3.2, we can assume that the $S^1\times \Bbb Z^k_p$-action
has a finite isotropy group $H$. Let $F$ be a component of $F(H,M)$
which contains a point $x$ whose isotropy group is $H$ (note that this
may not hold for every component of $F(H,M)$). We first claim that
$S^1\times
\Bbb Z_p^k$ preserves $F$. Since the $S^1$-subgroup always preserves
$F$, it suffices to show that $\Bbb Z_p^k$ preserves $F$. Otherwise,
since
$\Bbb Z_p^k$ preserves $F(H,M)$, then $F(H,M)$ must contain at least $p$
components.
Following the same arguments as in the proof of Lemma 2.1, one obtains
a contradiction.

Put $S^1\times \Bbb Z_p^r=(S^1\times \Bbb Z_p^k)/H$. Since $x\in F$,
$S^1\times \Bbb Z_p^r$ acts effectively on $F$. Note that $F$ is a
closed totally
geodesic submanifold of odd-dimension $<2n+1$. If $r=0$, then the proof
is
complete. If $r>0$, then we can apply the inductive assumption to
$(F,S^1\times \Bbb Z_p^r)$ and conclude that $S^1\times \Bbb Z_p^r$
has a subgroup $\Bbb Z_p^r$ whose fixed point set in $F$ is not empty.
This implies that $H\cong \Bbb Z_p^{k-r}$ and thus the $S^1\times \Bbb
Z_p^k$
has a $\Bbb Z_p^k$-subgroup whose fixed point set is not empty.

Finally, if $F_0$ is a $\Bbb Z^k_p$-fixed point component, then $\Bbb
Z_p^k$
acts effectively on a normal sphere space to $F_0$. If $\dim
(F_0)=2r+1$, then
from the proof of Lemma 2.2 $k\le n-r$ and thus $k\le n$. If $k=n$, then
$r=0$ and thus $\dim(F_0)=1$.
\qed\enddemo

\vskip2mm

\proclaim{Lemma 3.3} Let the assumptions be as in Theorem B with $k=n$.
Then there are closed totally geodesic submanifolds, $M^3\subset M^5
\subset\cdots \subset M^{2n+1}=M$ such that $M^{2i+1}$ admits an
isometric $\Bbb Z_p^i$-action.
\endproclaim

\vskip2mm

\demo{Proof} By Lemma 3.1, we can assume a subgroup $H\cong \Bbb Z^n_p$
whose fixed point set is a circle. Consider the isotropy representation
of $H$ on the normal space of the circle. Note that the normal space
has dimension $2n$. By now the rest of the proof follows exactly the
same argument as in the proof of Lemma 2.2.
\qed\enddemo

\vskip2mm

\demo{Proof of Theorem B}

By Theorems 1.2, it suffices to show that $M$ is a homotopy sphere
(note that for $n=1$, one may use [Ha]).

By Lemma 3.3, we obtain a sequence closed orientable totally
geodesic submanifolds,
$$M^3\subset M^5\subset\cdots \subset M^{2n+1}.$$
By Theorem 1.11, we conclude in order that $M^{2n-1},...,M^3$ are
simply connected.
Note that $M^3$ is a homotopy sphere. By Theorem 1.11 again, we then
conclude
that $M^5$ is $2$-connected and thus a homotopy sphere (by the
Poincar\'e duality and the Hurewicz theorem). By iterating the same
argument to $M^7,..., M^{2n+1}$,
we eventually conclude that $M^{2n+1}$ is homotopic to a sphere.
\qed\enddemo

\vskip4mm

\head 4. Proofs of Theorems C and D
\endhead

\vskip4mm

First, we need some preparation.

\vskip2mm

\proclaim{Lemma 4.1} Let $M$ be a closed simply connected 
$2n$-manifold of positive sectional curvature. If there is 
a closed submanifold $N$ of dimension at least $(n+1)$ 
such that $N$ is homotopic equivalent to 
a sphere or a complex projective 
space and the inclusion $i: N\hookrightarrow M^{2n}$ is at least 
an $n$-equivalence, then $M$ is homotopy equivalent to a sphere or a 
complex projective space.
\endproclaim

\vskip2mm

\demo{Proof} Observe that if $N$ is homotopic to a sphere, then 
by the Poincar\'e duality $M$ is a homotopy sphere. Hence, we 
can assume that $N$ is homotopic to a complex projective space. 
Since $\pi_1(M)=1$, it suffices to show that the cohomology
ring $H^*(M)$ is isomorphic to the truncated polynomial ring
$\Bbb Z[x]/(x^{n+1}=0)$ for a generator $x\in H^2(M)$.

Let $x\in H^2(M)$ be a generator. Note that $\pi_n(N)=0$. 
It is easy to see that $H^*(M)$ is isomorphic to $\Bbb Z$ 
up to degree $n$.

If $n$ is even, by the Poincar\'e duality $(x^{\frac n2})^2[M]=\pm 1$ 
and thus $H^*(M)$ is isomorphic to $\Bbb Z[x]/(x^{n+1}=0)$. 

If $n$ is odd,  write $x^ {\frac {n+1}2}=dy$ for some integer $d$,
where $y\in H^{n+1}(M)\cong \Bbb Z$ is a generator. Since $i^*(x)$ is
a generator of $H^2(N)$, the above implies that $di^*(y)=
(i^*(x))^ {\frac {n+1}2}$ is a generator of $H^{n+1}(N)\cong \Bbb Z$
and thus $d=\pm 1$. Consequently, $H^*(M)$ is also isomorphic to 
$\Bbb Z[x]/(x^{n+1}=0)$. 
\qed\enddemo

\vskip2mm

The following lemma was pointed to us by Burkhard Wilking.

\vskip2mm

\proclaim{Lemma 4.2} Let $M$ be closed simply connected 
$m$-manifold of positive sectional curvature. Assume that
 
\noindent (4.2.1) If $m=2n+1$, there is a closed totally
geodesic submanifold $N$ of codimension $2$ or 

\noindent (4.2.2) If $m=2n$, there are two closed totally geodesic 
submanifolds $N_1, N_2$ of codimenison $2$ such that
$\dim(N_1\cap N_2)=m-4>0$. 

Then $M$ is homotopy equivalent to a sphere or a complex projective space.
\endproclaim

\vskip2mm

Note that (4.2.1) also gives an alternative proof of Theorem B.

\vskip2mm

\demo{Proof} (4.2.1) By [Ha], we can assume $n>2$. By 
Theorem 1.2, it suffices to show that $M$ is a homotopy sphere. 
By (1.11.1) the inclusion $i: N\to M$ is $(2n-2)$-connected. 
In particular, $N$ is simply connected. By the Poincar\'e duality, 
$H^{1}(N)=H_{2n-2}(N)=0$ and thus $H^3(M)=H_{2n-2}(M)=0$. Thus
$H^{2n-4}(N)=H_3(N)=0$ if $m\ge 3$. This implies that $H^5(M)=H_{2n-4}(M)=0$.

By repeating this process one conclude that both $M$ and $N$ are 
homotopy spheres.

(4.2.2) Set $N=N_1\cap N_2$. Then $\dim(N)=2n-4>0$. By 
Theorem 1.11, $N_1\hookrightarrow M$ is $(2(n-2)+1)$-connected. 
By Lemma 4.1, it suffices to show that the cohomology ring $H^*(N_1)$ 
is isomorphic to $H^*(S^{2n-2})$ or $\Bbb Z[x]/(x^n=0)$.

Since $N$ is a closed geodesic submanifold of $N_i$ of codimension $2$,  
by Theorem 1.9 $N$ is connected. We now apply Theorem 1.11 
and conclude that the inclusion $i:
N\hookrightarrow  N_1$ is $(2n-4)$-connected. Let $[N]$ (resp. $[N_1]$) denote
the fundamental class of $N$ (resp. $N_1$). Since $\dim(N)=2n-4$,
it is clear that $i_*([N])$ generates $H_{2n-4}(N_1)$.

Let $x=PD(i_*([N]))\in H^2(N_1)$ be the Poincar\'e dual of $i_*([N])$.
It is standard to verify the composition homomorphism
$$H^j(N_1)@>i^* >>H^j(N)@>\cap [N] >> H_{2n-4-j}(N)@>i_*>>H_{2n-4-j}(N_1)
@>\cap [N_1]>> H^{j+2}(N_1)$$
is equal to $\cup x$, where $\cap$ (resp. $\cup$) is the cap (resp. cup)
product. Recall that $\cap [N]$ (resp. $\cap [N_1]$) is the Poincar\'e
dual map of $N$ (resp. $N_1$.)

Since $i: N\to N_1$ is a $(2n-4)$-equivalence, 
the homomorphism $H^j(N_1)@>i^* >>H^j(N)$ is an isomorphism for $j<2n-4$ and
an injection for $j=2n-4$. Similarly, $H_{2m-4-j}(N)@>i_*>>H_{2n-4-j}(N_1)$
is an isomorphism for $j>0$ and surjection for $j=0$. Therefore
the above homomorphism $\cup x: H^j(N_1)\to H^{j+2}(N_1)$ is an isomorphism
for $0<j<2m-4$. Therefore $H^{2m-1}(N_1)=0$ for all $1\le 2n-1\le 2n-3$.
By the universal coefficients theorem $H_*(N_1)$ is torsion free.

If $x=0$, then $\cup x=0$ and the above shows easily that $N_1$ is a homotopy
sphere.

If $x\ne 0$, by the above $H^2(N_1)\cong H_2(N_1)\cong \Bbb Z$ is generated
by
$x$. The above isomorphism shows that $x^{n-2}$ generates $H^{2n-4}(N_1)$.
By
Poincar\'e duality there is an element $y\in H^2(N_1)$ (since $\text{dim
}N_1=2n-2$) such that $x^{n-2}\cup y [N_1]=1$. Clearly, we may write $y=dx$
for
some integer $d$. Therefore $d=\pm 1$ and $x^{n-1}[N_1]=\pm 1$. This shows
that
the cohomology ring $H^*(N_1)$ is isomorphic to 
$\Bbb Z[x]/(x^n=0)$. 
\qed\enddemo

\vskip2mm

As a consequence of Proposition 2.3 and (4.2.2), we 
conclude a weak version of Problem 0.2.

\vskip2mm

\proclaim{Lemma 4.3} Let $M$ be closed simply connected 
$2n$-manifold of positive sectional curvature ($n\ge 2$). 
Assume $M$ admits an isometric $\Bbb Z_p\oplus \Bbb Z_p$-action 
($p\ge p(n))$ with fixed point set codimension $4$. Then $M$ 
is homotopic to a sphere or a complex projective space.
\endproclaim

\vskip2mm

\demo{Proof} Note that for $n=2$, Lemma 4.3 is included in 
Proposition 2.3. Hence, we can assume that $n\ge 3$.

Let $N\subseteq \text{Fix}(M,\Bbb Z_p\oplus \Bbb Z_p)$ be  
a fixed point component of codimension $4$. From the isotropy
representation, it is clear that there are two $\Bbb Z_p$ 
subgroups, $H_1$ and $H_2$, whose fixed point components $N_i$ 
contain $N$ and are both of codimension $2$. By now Lemma 
4.3 follows from (4.2.2). 
\qed\enddemo

\vskip2mm

The following is a consequence of Lemma 2.1 and Lemma 3.1.

\vskip2mm

\proclaim{Corollary 4.4} Let $N^n$ be a closed manifold of positive 
sectional curvature. 

\noindent (4.4.1) If $n$ is even and $\Bbb Z_p^k$ ($p\ge p(n)$) acts
isometrically on $N$, then $n\ge 2k$.

\noindent (4.4.2) If $n$ is odd and $S^1\times \Bbb Z_p^k$ ($p\ge p(n)$)
acts isometrically on $N$, then $n\ge 2k+1$. 
\endproclaim

\vskip2mm

\demo{Proof of Theorem C}

Observe that without loss of the generality, we may assume 
$k=[\frac{3n}4]+1$ or $2$. By Theorem 1.2 it suffices
to show $M$ is homotopic to a sphere or a complex projective 
space.

From the isotropy representation at a $\Bbb Z_p^k$-fixed point
(see Lemma 2.1), it is easy to see that the $\Bbb Z_p^k$-action has
an isotropy group $H\cong \Bbb Z_p$. Hence, we can assume an $H$-fixed 
point component $N$ such that on which $\Bbb Z_p^{k-1}=\Bbb Z_p^k/H$ 
acts effectively. By (4.4.1), $2(k-1)\le \text{dim}(N)\le \dim(M)-2$.

We shall proceed by induction on $n$. It turns out that in order to 
apply the inductive assumption (see (4.6)), we need to check the first 
four cases $7\le n\le 10$. 

Case 1. Assume $\dim(M)=14$ and thus $k=[\frac{3\cdot 7}4]+1=6$. Then 
$\dim(N)=10$ or $12$ and $\Bbb Z_p^5\cong \Bbb Z_p^6/H$ acts 
effectively on $N$. If $\dim(N)=10$, then by Theorem A, 
$N$ is homeomorphic to $S^{10}$ or $\Bbb CP^5$. By Theorem 1.11, 
$N\hookrightarrow M$ is $7$-connected. Applying Lemma 4.1 to 
$(M,N)$, we conclude that $M$ is homotopic to $S^{14}$ or $\Bbb CP^7$.

If $\dim(N)=12$, by Theorem 1.11 $N\hookrightarrow M$ is $11$-connected.
By Lemma 4.1 it suffices to prove that $N$ is homotopic to $S^{12}$ or
$\Bbb CP^6$.

Consider $(N,\Bbb Z_p^5)$. First, $\Bbb Z_p^5$ has a subgroup $H_1$ 
with a fixed point component $N_1$ of dimension equals to $8$ or $10$ and 
$\Bbb Z_p^4\cong \Bbb Z_p^5/H_1$ acts effectively on $N_1$. If $\dim(N_1)=8$, then by Theorem, A $N_1$ is homeomorphic to $S^8$ or $\Bbb CP^4$. By 
Theorem 1.11, $N_1\hookrightarrow N$ is $5$-connected. By Lemma 4.1, 
we then conclude that $N$ is homotopic to $S^{12}$ or $\Bbb CP^6$.

If $\dim(N_1)=10$, then we are in a situation that $(M,\Bbb Z^{6})$ 
has a $\Bbb Z_p\oplus \Bbb Z_p$-subgroup whose fixed point set has codimension
$4$. 
We then conclude the desired result from Lemma 4.3.

Case 2. Assume $\dim(M)=16$ and thus $k=7$. Then $\dim(N)=12$ or $14$. 
If $\dim(N)=12$, then by Theorem A, $N$ is homeomorphic to $S^{12}$ 
or $\Bbb CP^6$. By Theorem 1.11, $N\hookrightarrow M$ is $9$-connected. 
Now Lemma 4.1 applies to conclude that $M$ is homotopic to $S^{16}$ 
or $\Bbb CP^8$.

If $\dim(N)=14$, then by Case 1 we then conclude that $N$ is homotopic 
to $S^{14}$ or $\Bbb CP^7$. By Theorem 1.11 and Lemma 4.1, we then conclude 
that $M$ is homotopic to $S^{16}$ or $\Bbb CP^8$.

Case 3. Assume $\dim(M)=18$ and thus $k=[\frac{3\cdot 9}4]+2=8$. Then
$\dim(N)=14$ or $16$. 
If $\dim(N)=14$, then by Theorem A $N$ is homeomorphic to $S^{14}$ 
or $\Bbb CP^7$. By Theorem 1.11, $N\hookrightarrow M$ is $11$-connected. 
Now Lemma 4.1 applies to conclude that $M$ is homotopic to $S^{18}$ 
or $\Bbb CP^9$.

If $\dim(N)=16$, then by Case 2 we then conclude that $N$ is homotopic 
to $S^{16}$ or $\Bbb CP^8$. By (4.2.1) and Theorem 1.11, we then conclude 
that $M$ is homotopic to $S^{18}$ or $\Bbb CP^9$.

Case 4. Assume $\dim(M)=20$ and thus $k=[\frac{3\cdot 10}4]+2=9$. 
Then $\dim(N)=16$ or $18$. 
If $\dim(N)=16$, then by Theorem A, $N$ is homeomorphic to $S^{14}$ 
or $\Bbb CP^7$. By Theorem 1.11, $N\hookrightarrow M$ is $13$-connected. 
Now Lemma 4.1 applies to conclude that $M$ is homotopic to $S^{20}$ 
or $\Bbb CP^{10}$.

If $\dim(N)=18$, then by Case 3 we then conclude that $N$ is homotopic 
to $S^{16}$ or $\Bbb CP^8$. By Theorem 1.11 and Lemma 4.1, we then conclude 
that $M$ is homotopic to $S^{20}$ or $\Bbb CP^{10}$.

We now proceed the rest of the proof by induction on $\dim(M)=2n$ 
starting with $n=10$. Assume that Theorem C holds for $n\ge 10$. 

Consider $\dim(M)=2(n+1)$. 

\noindent a. Assume $n\ne 0\, \,(\text{mod}\,\,4)$ (i.e., $n+1\ne 1(\text{mod
}
\,4)$). Then $k=[\frac{3(n+1)}4]+1=[\frac{3n}4]+2$.
As in the above, there is a subgroup $H$ isomorphic to $\Bbb Z_p$ with a
fixed
point component $N$ satisfying 
$$14\le 2\left(\left[\frac{3n}4\right]+1\right)\le \dim(N)\le 2n \qquad ( 
n\ge 9)\tag 4.5$$
since $\Bbb Z_p^{[\frac{3n}4]+1}\cong \Bbb Z^{[\frac{3n}4]+2}/H$ 
acts effectively on $N$. 

If $n\ne 1, 2 \,\, (\text{mod } 4)$, by (4.5) we can apply the 
inductive assumption to $(N,\Bbb Z_p^{[\frac{3n}4]+1})$ and 
conclude that $N$ is homotopic to a sphere or 
a complex projective space. Since $\dim(N)\ge 2[\frac {3n}4]+2\ge n+2$, 
$[2\dim(N)-2(n+1)+1]\ge n+3$. Hence, by Theorem 1.11 $N\hookrightarrow M$ is at least $(n+3)$-connected. By Lemma 4.1 we conclude that $M$
is homotopic to $S^{2(n+1)}$ or $\Bbb CP^{n+1}$.

If $n=1 \text{ or }2 \,\, (\text{mod } 4)$, then $(N,\Bbb
Z_p^{[\frac{3n}4]+1})$ has 
an isotropy group $H_1$ isomorphic to $\Bbb Z_p$ with a fixed point 
component $N_1$ satisfying 
$$14\le 2\left[\frac{3n}4\right]\le \dim(N_1)\le \dim(N)-2\le 2(n-1)\qquad
(n\ge 10)\tag 4.6$$
If $\dim(N_1)\le 2(n-2)$, then $[\frac{3n}4]\ge [\frac{3(n-2)}4]+1$. By 
(4.6) we can apply the inductive assumption to $(N_1,\Bbb Z_p^{[\frac{3n}4]})$ to conclude that $N_1$ is homotopic to a sphere or 
a complex projective space, and therefore by Lemma 4.1,
$N$ is homotopic to a sphere or a complex projective space. Similarly,
$M$ is homotopic to $S^{2(n+1)}$ or $\Bbb CP^{n+1}$. 

If $\dim(N_1)=2(n-1)$, then we are in a situation that $(M,\Bbb Z_p^{[\frac 
{3n}4]+2})$ has a $\Bbb Z_p\oplus \Bbb Z_p$ subgroup of fixed point 
set codimension $4$. We then conclude the desired result by applying Lemma
4.2. 

\noindent b. Assume $n=0\, \,(\text{mod}\,\,4)$. Then
$[\frac{3(n+1)}4]+1=\frac{3n}4+1$.
Since $n\ne 1, 2\, \,(\text{mod}\,\,4)$, we can 
apply the inductive assumption to conclude that $N$ is homotopic to a sphere or a complex projective space. Since $\dim(N)\ge 2([\frac {3n}4]+1)\ge n+4$ 
(because $n\ge 10$), $[2\dim(N)-2(n+1)+1]\ge n+3$. 
Hence, by Theorem 1.11 $N\hookrightarrow M$ is at least $(n+3)$-connected. 
Once again by Lemma 4.1, we conclude that $M$ is homotopic to 
$S^{2(n+1)}$ or $\Bbb CP^{n+1}$.
\qed\enddemo

\vskip2mm

\demo{Proof of Theorem D}

We shall follow the strategy in the proof of Theorem C. 
By Theorem 1.2, it suffices to show that $M$ is a homotopy sphere.
 
We first observe that by Sublemma 3.2, there is 
an $H$-fixed point component $N$ with $H\cong \Bbb Z_p$ and 
$(S^1\times \Bbb Z_p^k)/H$ acts effectively on $N$. By (4.4.2), 
$2(k-1)+1\le \text{dim}(N)\le \dim(M)-2$.
 
Case 1. Assume $\dim(M)=2\cdot 7+1=15$ and thus $k=[\frac{3\cdot 7}4]+1=6$. 
From the above, $\dim(N)=11$ or $13$ and $S^1\times \Bbb Z_p^6/H$ acts 
effectively on $N$. If $\dim(N)=13$, then by (4.2.1) 
$M$ is homeomorphic to $S^{15}$. If $\dim(N)=11$, then by Theorem B, $N$ is 
homeomorphic to $S^{11}$ (note that in this case, $(S^1\times \Bbb Z_p^6
)/H\cong S^1\times \Bbb Z_p^5$). By Theorem 1.11, $N\hookrightarrow M$ is 
$8$-connected. Hence, by the Poincar\'e duality $M$ is homotopic to
$S^{15}$.

Case 2. Assume $\dim(M)=17$ and $k=7$. Then $\dim(N)=13$ 
or $15$ and $(S^1\times \Bbb Z_p^6)/H$ acts 
effectively on $N$. If $\dim(N)=15$, then by (4.2.1) $M$ is 
homeomorphic to $S^{17}$. If $\dim(N)=13$, then by Theorem B $N$ is 
homeomorphic to $S^{13}$. By Theorem 1.11, $N\hookrightarrow M$ is 
$10$-connected. By the Poincar\'e duality $M$ is homotopic to $S^{17}$.

We now proceed the rest of the proof by induction on 
$\dim(M)=2n+1$ starting with $n=8$. 

Consider $\dim(M)=2(n+1)+1$ ($n\ge 8$).

\noindent a. Assume $n+1=1\,(\text{mod } 4)$ and thus $k=[\frac{3(n+1)}4]+2$.

Note that $(S^1\times \Bbb Z_p^{[\frac{3(n+1)}4]+2})/H$ contains a subgroup 
isomorphic to $S^1\times \Bbb Z_p^{[\frac{3(n+1)}4]+1}$. 
Since $n\ge 8$, 
$$15\le 2\left(\left [\frac{3(n+1)}4\right]+1\right)+1\le \dim(N)\le 2n+1.\tag
4.7$$
If $\dim(N)=2n+1$, then by (4.2.1) $M$ is homeomorphic to 
a sphere. Note that 
$$\left[\frac{3(n+1)}4\right]+1\ge \left[\frac{3(n-1)}4\right]+2.$$
If $\dim(N)\le 2(n-1)+1$, then by (4.7) and Case 1 we can apply the 
inductive assumption 
and conclude that $N$ is homotopic to a sphere. Since 
$$2\dim(N)-(2(n+1)+1)+1 \ge 
4\left[\frac{3(n+1)}4\right]+2-(2(n+1)+1)+1=n+3,$$
by Theorem 1.11, $N\hookrightarrow M$ is at least $(n+3)$-connected.
Hence, by the Poincar\'e inequality, $M$ is also homotopic to a sphere.

\noindent b. Assume $n+1\ne 1 \,\, (\text{mod } 4)$ and thus 
$k=[\frac{3(n+1)}4]+1$. Since $n\ge 9$ (note that $n=8$ is considered in
a), 
$$15\le 2\left(\left [\frac{3(n+1)}4\right]+1\right)+1\le \dim(N)\le 2n+1, \tag
4.8$$
and $(S^1\times \Bbb Z_p^{[\frac{3(n+1)}4]+1})/H$ acts effectively on 
$N$. If $\dim(N)=2n+1$, then by (4.2.1) $M$ is homeomorphic to 
a sphere. Note that 
$$\left[\frac{3(n+1)}4\right]\ge \cases [\frac{3(n-1)}4]+1 & n-1\ne 1
(\text{mod } 4)\\
[\frac{3(n-1)}4]+2 & n-1 =1 (\text{mod } 4)\endcases.$$
If $\dim(N)\le 2(n-1)+1$, then by (4.8) and Case 1 we can apply the inductive assumption and conclude that $N$ 
is homeomorphic to a sphere. Since 
$$2\dim(N)-(2(n+1)+1)+1\ge 
4\left[\frac{3(n+1)}4\right]+2-(2(n+1)+1)+1=n+3,$$
by Theorem 1.11 $N\hookrightarrow M$ is at least $(n+3)$-connected.
Hence, by the Poincar\'e duality $M$ is also homotopic to a sphere.
\qed\enddemo

\enddocument